\documentclass[12pt,a4paper]{article}
\usepackage{amsmath,amssymb,amsthm,mathptmx}

\author{Anastasia V. Parusnikova}
\title{4-Dimensional Power Geometry and its Application to $P_1$ -- $P_5$}

\date{}

\begin{document}
\maketitle

\begin{abstract}
    %
In the first section of this work we introduce 4-dimensional Power Geometry for second-order ODEs of a polynomial form. In the next five sections we apply this construction to the first five Painlev\'{e} equations. 


\textbf{ MSC classes:} 	34m25, 34m55.


    %
\end{abstract}

\section{4-dimensional Power Geometry}
We introduce definitions and notations of 4-dimensional Power Geometry analogous to the way it has been done in two and three-dimensional cases in \cite{BrUmn} and \cite{Br-3d}.

Let us be given a second order differential equation of the form
\begin{equation}
 \label{eq1.1}
f(z,w,w',w'')=0,
\end{equation}
where $f(z,w,w',w'')$~is a polynomail, $z$ is an independent, $w$ is a dependent variable, $w'=dw/dz$.

To each monomial $a(z,w)$ in polynomial (\ref{eq1.1}) we put in correspondance its \textit{four-dimensional exponent}  $\mathbf{Q}(a(z,w))=(q_1,q_2,q_3,q_4)$ according to the following rule:
$$
\mathbf{Q}(c z^r w^s)=(r,s,0,0);~
\mathbf{Q}\left(\frac{d w}{d z}\right)=(0,0,1,0);~
\mathbf{Q}\left(\frac{d^2 w}{d z^2}\right)=(0,0,0,1);
$$
$$\mathbf{Q}(a(z,w)b(z,w))=\mathbf{Q}(a(z,w))+\mathbf{Q}(b(z,q)).$$
The set of all exponents of the monomials in polynomial $f(z,w)$ is called
\textit{a support of a differential sum} $f(z,w)$ and is denoted as
$\mathbf{\tilde{S}}(f)$. The convex hull ${\Gamma}(f)$
of the support $\mathbf{\tilde{S}}(f)$ is called \textit{polyhedron of a differential sum} $f(z,w)$, the boundary
 $\partial{\Gamma}(f)$ consists of the vertices
${\Gamma}^{(0)}_j$, edges ${\Gamma}^{(1)}_j$,
two-dimensional faces ${\Gamma}^{(2)}_j$ and three-dimensional faces ${\Gamma}^{(3)}_j$ (hyper faces).

We put in correspondence \textit{a truncated equation} $f^{(d)}_j(z,w,w',w'')=0$  to every edge ${\Gamma}^{(d)}_j$, where $f^{(d)}_j(z,w,w',w'')=\sum\limits_{\displaystyle{a_s(z,w):\mathbf{Q}(a_s(z,w))\in {\Gamma}^{(d)}_j}}a_s(z,w).$

Let us compare a four-dimensional construction with three and two-dimensional ones. 

In three-dimensional case we put in correspondance to each differential monomial $a(z,w)$ its \textit{three-dimensional exponent}  $\mathbf{Q}(a(z,w))=(q_1,q_2,q_3)$ according to the following rule:
$$
\mathbf{Q}(c z^r w^s)=(r,s,0);~
\mathbf{Q}\left(\frac{d^l w}{d z^l}\right)=(0,1,l);
$$
$$\mathbf{Q}(a(z,w)b(z,w))=\mathbf{Q}(a(z,w))+\mathbf{Q}(b(z,w)).$$

In three-dimensional case we put in correspondance to each differential monomial $a(z,w)$ its \textit{three-dimensional exponent}  ${Q}(a(z,w))=(q_1,q_2,q_3)$ according to the following rule:
$$
{Q}(c z^r w^s)=(r,s);~
{Q}\left(\frac{d^l w}{d z^l}\right)=(-l,1);
$$
$${Q}(a(z,w)b(z,w))={Q}(a(z,w))+{Q}(b(z,w)).$$

Let us redefine a notion of \textit{order of the function}. Let us be given a function $\psi(z)$, for which infinity is not an accumulation point of poles and zeroes, $z=re^{i \varphi}$. We call $$p_+(\psi(z),\varphi)=\displaystyle{{\lim\limits_{r \rightarrow \infty }}\frac{\ln|\psi(re^{i \varphi})|}{\ln|r|}}$$ an order of a function on a ray with direction $\varphi$  ($z \rightarrow \infty$) if this limit exists.

Let us be given a function $\psi(z)$, for which zero is not an accumulation point of poles and zeroes of the function, $z=re^{i \varphi}$. We call $$p_-(\psi(z),\varphi)=\displaystyle{{\lim\limits_{r \rightarrow 0 }}\frac{\ln|\psi(re^{i \varphi})|}{\ln|r|}}$$  an order of a function on a ray with direction $\varphi$  ($z \rightarrow 0$) if this limit exists.

These orders of a function $f(z)$ are of special interest if they coincide for $\varphi\in(\varphi_1,\varphi_2)$, i.e. for points in some sector on Riemann surface of the logarithm.


Thus we see that in two-dimensional case construction of the support assumes automatically that formal asymptotics of the solutions to the equation satisfy the following condition:  order of the derivative of the function is one less that an order of the function (for example, the function $z^n$ satisfies this condition). 

In the three-dimensional structure is assumed that at each differentiation order of the function is changed by $\gamma_1 \in \mathbb{R}$. Functions with $\gamma_1\neq 1$ as in 2-D variant also exist: consider $p_-(\sin z, 0)$ and $p_-(\cos z, 0)$).

Introduction of the fourth coordinate differential monomials exponents included in the second order ODEs
is due to the following condition on the asymptotic behavior of the possible solutions: the first  differentiation changes order of a function by $\gamma_1 \in \mathbb{R} $, the second differentiation changes order of a function by $\gamma_2 \in \mathbb {R}$.
In 2-D construction $\gamma_1=\gamma_2=1$, in 3-D construction $\gamma_1=\gamma_2$. Functions with $\gamma_1\neq\gamma_2$ exist: consider $p_-(z+z^{\alpha},0)$, $p_-(1+z^{\alpha-1},0)$ and $p_-(z^{\alpha-2},0)$, $\alpha >1$.

We formulate a necessary condition of the fact that a truncated equation corresponding to hyper-face in four-dimensional case can have a solution which is an asymptotic form of a solution to the initial equation.

\textbf{Assertion 1.} 
Let us be given a differential equation of the second order. We consider a truncated equation corresponding to a hyper face ${\Gamma}^{(3)}$ with an external normak $\mathbf{N}=(n_1,~n_2,~n_3,~n_4)$. If $n_1 = 0$, then a truncated equation has no solution of finite order which can be an asymptotic form of the solution to the initial equation.

\textbf{Proof.} Negative proof. Let the hyper face have an external normal $N=(0,~n_2,~n_3,~n_4)$, hyperplane containing the face has an equation
 \begin{equation}
                 \label{plane}                                                                                                 
 n_2 q_2+ n_3 q_3+ n_4 q_4+e=0.                          
                 \end{equation}
 There exist 4 points 
${Q}_i=(q_{1i},~q_{2i},~q_{3i},~q_{4i}),~i=1,\ldots,~4$ in the face, which lie on the hyperplane (\ref{plane}) but do not lie in a plane of smaller dimension. That means that the rank of the matrix $Q=(q_{ij}),~i,~j=1,~\ldots,~4$ is equal to three.

All the points $Q_i$ satisfy an equation (\ref{plane}). 
We subtract from each equation the first equation and obtain an equation with the matrix
$$\tilde{Q}=(q_{ij}-q_{i1}),~i=1,\ldots,4,~j=1,~2,~3.$$  As $\tilde{Q} N^T=0$, and $N=(0,~n_2,~n_3,~n_4)$,  we obtain that the matrix $\hat{Q}=(q_{ij}-q_{i1}),~i=2,\ldots,4,~j=1,~2,~3$ has a rank less than three. If we assume that it is less than two we arrive at a contradiction, we obtain that the rank of $\hat{Q} $ is equal to two , i.e. its columns are linearly dependent.

If a truncated equation corresponding to a hyper face has a solution an order of which is equal to $\gamma$, an order of the first derivative is equal to $\gamma_1$, an order of the second derivative is equal to  $\gamma_2$, then the points ${Q}_i,~i=1,\ldots,~4$ satisfy the system $$q_{1i}+q_{2i}\gamma+q_{3i}\gamma_1+q_{4i} \gamma_2+f=0,~i=1,\ldots,~4.$$ We subtract from the second, third and fourth equation the first one and obtain a system $$\hat{Q}x=y,$$ where $x=(\gamma,~\gamma_1,~\gamma_2)^T,$ $y=(q_{11}-q_{12},~q_{11}-q_{13},~q_{11}-q_{14})^T$. But the rank of the matrix $\hat{Q}$ is equal to 2, and a column $y$ is a linear combination of the columns of the matrix $\hat{Q}$. We obtain that rank of the matrix $\tilde{Q}$ is equal to 2: this leads us to contradiction.

\section{The first Painlev\'{e} equation}
The first Painlev\'{e} equation
\begin{equation*}
w''=6 w^2+z
\end{equation*}
has one singular point $z=\infty$.

All three points of the support, of course, lie in the plane of dimension 3, so the use of the methods of four-dimensional power geometry does not make sense.

\section{The second Painlev\'{e} equation}
We consider the second Painlev\'{e} equation
\begin{equation*}
w''=2w^3+zw+\alpha,
\end{equation*}
where $\alpha$ is a complex parameter.
The equation  has one singular point $z=\infty$.

The support of the equation 
consists of four points: $(0,0,0,1)$, $(0,3,0,0)$ $(1,1,0,0)$ and $(0,0,0,0)$, all these points lie in a hyperplane $q_3=0$, an external normal to it is equal either to $(0,0,1,0)$ or to $(0,0,-1,0)$. According to the assertion 1 we obtain that the truncated solution has no solutions the leading term of which with its derivatives have a finite order.

\section{The third Painlev\'{e} equation $\alpha\beta\gamma\delta \neq 0$}
We consider the third Painlev\'{e} equation 
\begin{equation}
\label{P3} 
w''=\frac{(w')^2}{w}-\frac{w'}{z}+\frac{\alpha
w^2+ \beta}{z}+ \gamma w^3+\frac{\delta}{w},
\end{equation}
where $\alpha, \beta, \gamma, \delta$
are complex paramters, the equation (\ref{P3}) has two singular points: $z=0$ and~$z=\infty$.

We rewrite the equation (\ref{P3}) in a form of a differential sum:
\begin{equation}
\label{P3_mod}
 f_3(z,w)\stackrel{def}{=}-z w  w''+z \left(w'\right)^2
- w w'
+\alpha w^3+\beta w+\gamma z w^4+\delta z=0.
\end{equation}

The four-dimensional support of the equation (\ref{P3_mod}) consists of the points
\begin{equation*}
 \mathbf{Q}_1=(1,0,0,0),~\mathbf{Q}_2=(1,4,0,0),~\mathbf{Q}_3=(0,1,0,0),~\mathbf{Q}_4=(0,3,0,0),$$ $$
\mathbf{Q}_5=(0,1,1,0),~\mathbf{Q}_6=(1,0,2,0),~\mathbf{Q}_7=(1,1,0,1).
\end{equation*}
A convex hull of the equation (\ref{P3_mod}) is a polygon with vertices
$\mathbf{Q}_1,\ldots,~\mathbf{Q}_7$.

 Its hyper faces are
$
 {\Gamma}_1^{(3)}=\textrm{conv}(\mathbf{Q}_1,\ldots,\mathbf{Q}_{6}),~
{\Gamma}_2^{(3)}=\textrm{conv}(\mathbf{Q}_1,\mathbf{Q}_{3},\mathbf{Q}_{5},\mathbf{Q}_{6},\mathbf{Q}_{7}),$
$ {\Gamma}_3^{(3)}=\textrm{conv}(\mathbf{Q}_1, \mathbf{Q}_{2},\mathbf{Q}_{6},\mathbf{Q}_{7},)$,
${\Gamma}_4^{(3)}=\textrm{conv}(\mathbf{Q}_1,\mathbf{Q}_2,\mathbf{Q}_3,\mathbf{Q}_{4},
\mathbf{Q}_{7}),$\newline
${\Gamma}_5^{(3)}=\textrm{conv}(\mathbf{Q}_2,\mathbf{Q}_{4},\ldots,\mathbf{Q}_{7}),$
${\Gamma}_6^{(3)}=\textrm{conv}(\mathbf{Q}_3,
\mathbf{Q}_{4},\mathbf{Q}_{5},\mathbf{Q}_{7}).$

According to the Assertion 1 we consider only the following 3D faces (i.e. the faces the first coordinate of the external normal to which is not equal to zero):

1. $
{\Gamma}_2^{(3)}$ lies in the plane $q_1+q_2-q_4-1=0$ with an external normal $N_2=(-1,-1,0,1)$.
The truncated equation corresponding to the face
$$
 -z w  w''+z \left(w'\right)^2
- w w'
+\beta w+\delta z=0
$$
has already been considered as a truncated equation corresponding to a 2D face in 3D case and as a truncated equation corresponding to an edge in 2D case.

2. $
{\Gamma}_3^{(3)}$ lies in the plane $q_1=1$ with an external normal $N_3=(-1,0,0,0)$.
The truncated equation corresponding to the face
$$
 -z w  w''+z \left(w'\right)^2
+\gamma z w^4+\delta z=0
$$
has already been considered as a truncated equation corresponding to a 2D face in 3D case.

3.
$
{\Gamma}_5^{(3)}$ lies in the plane $q_1-q_2-2q_3-3q_4+3=0$ with an external normal $N_5=(-1,1,2,-3)$.
The truncated equation corresponding to the face
$$
 -z w  w''+z \left(w'\right)^2
- w w'
+\alpha w^3+\gamma z w^4=0
$$
has already been considered as a truncated equation corresponding to an edge in 2D case.

4. $
{\Gamma}_6^{(3)}$ lies in the plane $q_1-q_4=0$ with an external normal $N_6=(-1,0,0,1)$.
 The truncated equation corresponding to the face
$$
 -z w  w''
- w w'
+\alpha w^3+\beta w=0.
$$
This equation has not been considered neither in 2D nor in 3D case.

The external normals to the faces
${{\Gamma}}_1^{(3)}$ ($q_4=0$, ${N}_1=(0,0,0,-1)$) and
${{\Gamma}}_4^{(3)}$  ($q_3=0$, ${N}_4=(0,0,-1,0)$)  do not satisfy Assertion 1.

\section{The fourth Painlev\'{e} equation $\alpha\beta \neq 0$}
We consider the fourth Painlev\'{e} equation 
\begin{equation}
\label{P4} 
w''=\frac{(w')^2}{w}-\frac{3}{2}w^3+4zw^2+2(z^2-\alpha)w+\frac{\beta}{w},
\end{equation}
where $\alpha, \beta,$
are complex parameters.

We rewrite the equation (\ref{P4}) in a form of a differential sum:
\begin{equation}
\label{P4_mod}
 f_4(z,w)\stackrel{def}{=}-2 w  w''+ \left(w'\right)^2
+3w^4+8zy^4+4(z^2-\alpha)w^2+2\beta=0.
\end{equation}

The 4D support of the equation (\ref{P4_mod}) consists of the points
\begin{equation}
\label{points_4}
 \mathbf{Q}_1=(0,0,0,0),~\mathbf{Q}_2=(0,2,0,0),~\mathbf{Q}_3=(2,2,0,0),~\mathbf{Q}_4=(1,3,0,0),$$ $$
\mathbf{Q}_5=(0,4,0,0),~\mathbf{Q}_6=(0,0,2,0),~\mathbf{Q}_7=(0,1,0,1).
\end{equation}
A convex hull of the equation (\ref{P3_mod}) is a polygon with vertices
$\mathbf{Q}_1,~ \mathbf{Q}_3,~ \mathbf{Q}_5,~\mathbf{Q}_6,~\mathbf{Q}_7$ with hyper faces
$
 {\Gamma}_1^{(3)}=\textrm{conv}(\mathbf{Q}_1,\ldots,\mathbf{Q}_{6}),$
$
{\Gamma}_2^{(3)}=\textrm{conv}(\mathbf{Q}_3,\ldots,\mathbf{Q}_{7}),$\newline
$ {\Gamma}_3^{(3)}=\textrm{conv}(\mathbf{Q}_1, \mathbf{Q}_{3},\mathbf{Q}_{6},\mathbf{Q}_{7}),$ $
{\Gamma}_4^{(3)}=\textrm{conv}(\mathbf{Q}_1,\ldots,\mathbf{Q}_5,\mathbf{Q}_{7}),$ \newline$
{\Gamma}_5^{(3)}=\textrm{conv}(\mathbf{Q}_1,\mathbf{Q}_{2},\mathbf{Q}_5,\mathbf{Q}_6,\mathbf{Q}_{7}).
$

According to the Assertion 1 we consider only the following 3D faces:

1. $
{\Gamma}_2^{(3)}$ lies in the plane $q_1+q_2+2q_3+3q_4-4=0$ with an external normal $N_2=(1,1,2,3)$.
The truncated equation corresponding to the face
$$
-2 w  w''+ \left(w'\right)^2
+3w^4+8zy^4+4z^2w^2=0$$
has already been considered as a truncated equation corresponding to a 2D face in 3D case.

2.
$
{\Gamma}_3^{(3)}$ lies in the plane $q_1-q_2+q_4=0$ with an external normal $N_3=(1,-1,0,1)$.
The truncated equation corresponding to the face
$$
-2 w  w''+ \left(w'\right)^2
+4z^2w^2+2\beta=0$$
has already been considered as a truncated equation corresponding to a 2D face in 3D case.

3. $
{\Gamma}_5^{(3)}$ lies in the plane $q_1=0$ with an external normal $N_5=(-1,0,0,0)$.
The truncated equation corresponding to the face
$$
-2 w  w''+ \left(w'\right)^2
+3w^4-4\alpha w^2+2\beta=0$$
has already been considered as a truncated equation corresponding to a 2D face in 3D case.

The external normals to the faces
${{\Gamma}}_1^{(3)}$ with equation $q_4=0$ and
${{\Gamma}}_4^{(3)}$ with equation $q_3=0$ (${{N}}_1=(0,0,0,-1)$
and $N_4=(0,0,-1,0)$) do not satisfy the  above conditions.

\section{The fifth Painlev\'{e} equation (the case $\delta \neq 0$)}
We consider the fifth Painlev\'{e} equation
\begin{equation*}
 w''=\left( \frac{1}{2w}+\frac{1}{w-1}
\right)\left(w'\right)^2-\frac{w'}{z}+\frac{(w-1)^2}{z^2}\left(\alpha
w+ \frac{\beta}{w}\right)+ \frac{\gamma w}{z}+\frac{\delta w
(w+1)}{w-1},
\end{equation*}
where $\alpha, \beta, \gamma, \delta$
are complex parameters.
The equation  has two singular points $z=0$ and~$z=\infty$.

We represent the fifth Painlev\'{e} equation  in a form of a differential sum:
\begin{equation}
\label{P5_mod}
 f(z,w)\stackrel{def}{=}-z^2 w (w-1) w''+z^2 \left(
 \frac{3}{2}w-\frac{1}{2}\right)\left(w'\right)^2
-z w (w-1) w'+$$
$$ +(w-1)^3(\alpha w^2+\beta)+
\gamma z w^2 (w-1)+\delta z^2 w^2 (w+1)=0.
\end{equation}

The 4D support of the equation (\ref{P5_mod}) consists of the points
\begin{equation*}
 \mathbf{Q}_1=(2,3,0,1),~\mathbf{Q}_2=(2,2,0,1),~\mathbf{Q}_3=(2,3,2,0),~\mathbf{Q}_4=(2,2,2,0),$$ $$
\mathbf{Q}_5=(1,2,1,0),~\mathbf{Q}_6=(1,3,1,0),~\mathbf{Q}_7=(0,0,0,0),~\mathbf{Q}_{8}=(0,5,0,0),$$
 $$
\mathbf{Q}_{9}=(1,2,0,0),~\mathbf{Q}_{10}=(1,3,0,0),$$ $$\mathbf{Q}_{11}=(2,2,0,0),~\mathbf{Q}_{12}=(2,3,0,0), $$ $$\mathbf{Q}_{13}=(0,1,0,0),~\mathbf{Q}_{14}=(0,2,0,0),\mathbf{Q}_{15}=(0,3,0,0),~\mathbf{Q}_{16}=(0,4,0,0).
\end{equation*}
A convex hull of the equation (\ref{P5_mod}) is a polygon with vertices
$\mathbf{Q}_1,~\mathbf{Q}_2,~\mathbf{Q}_3,$ $\mathbf{Q}_{4},$ $\mathbf{Q}_7,$ $\mathbf{Q}_{8},$ $
\mathbf{Q}_{11},$ $\mathbf{Q}_{12}$ with hyper faces
\begin{equation*}
 {\Gamma}_1^{(3)}=\textrm{conv}(\mathbf{Q}_1,\mathbf{Q}_{2},\mathbf{Q}_{3},\mathbf{Q}_{4},\mathbf{Q}_{11},\mathbf{Q}_{12}),\,
{\Gamma}_2^{(3)}=\textrm{conv}(\mathbf{Q}_2,\mathbf{Q}_{4},\mathbf{Q}_{7},\mathbf{Q}_{11}),$$
$$ {\Gamma}_3^{(3)}=\textrm{conv}(\mathbf{Q}_1, \mathbf{Q}_{2},\mathbf{Q}_{7},\mathbf{Q}_{12},\mathbf{Q}_{13},\ldots,\mathbf{Q}_{16}),\,
{\Gamma}_4^{(3)}=\textrm{conv}(\mathbf{Q}_1,\mathbf{Q}_3,\mathbf{Q}_8,\ldots,
\mathbf{Q}_{12}),$$ $$
{\Gamma}_5^{(3)}=\textrm{conv}(\mathbf{Q}_3,\ldots,
\mathbf{Q}_{12},\mathbf{Q}_{13},\ldots,\mathbf{Q}_{16}),\,
{\Gamma}_6^{(3)}=\textrm{conv}(\mathbf{Q}_1,\ldots,
\mathbf{Q}_{8},\mathbf{Q}_{13},\ldots,\mathbf{Q}_{16}).
\end{equation*}

According to the Assertion 1 we consider only the following 3D faces:

1. $
{\Gamma}_1^{(3)}$ lies in the plane $q_1=2$ with an external normal $N_1=(1,0,0,0)$.
The truncated equation corresponding to the face
$$
-z^2 w (w-1) w''+z^2 \left(
 \frac{3}{2}w-\frac{1}{2}\right)\left(w'\right)^2
 +\delta z^2 w^2 (w+1)=0$$
has already been considered as a truncated equation corresponding to a 2D face in 3D case.

2. $
{\Gamma}_2^{(3)}$ lies in the plane $q_1-q_2-q_3-q_4=0$ with an external normal $N_2=(1,-1,-1,-1)$.
The truncated equation corresponding to the face
$$
z^2 w  w''-\frac{1}{2}z^2\left(w'\right)^2
-\beta+\delta z^2 w^2 =0$$
has already been considered as a truncated equation corresponding to a 2D face in 3D case.

3.
$
{\Gamma}_4^{(3)}$ lies in the plane $q_1+q_2+q_3+q_4-5=0$ with an external normal $N_4=(1,1,1,1)$.
The truncated equation corresponding to the face
$$
-z^2 w^2 w''+ \frac{3}{2}z^2w\left(w'\right)^2
+\alpha w^5+
\delta z^2 w^3 =0$$
has already been considered as a truncated equation corresponding to a 2D face in 3D case.

4. $
{\Gamma}_6^{(3)}$ lies in the plane $q_1-q_3-2q_4=0$ with an external normal $N_6=(1,0,-1,-2)$.
The truncated equation corresponding to the face
$$
-z^2 w (w-1) w''+z^2 \left(
 \frac{3}{2}w-\frac{1}{2}\right)\left(w'\right)^2
-z w (w-1) w'
 +(w-1)^3(\alpha w^2+\beta)=0$$
has already been considered as a truncated equation corresponding to a 2D face in 3D case, this equation can be solved directly.

The external normals to the faces
${{\Gamma}}_3^{(3)}$ with equation $q_3=0$ and
${{\Gamma}}_5^{(3)}$ with equation $q_4=0$ (${N}_3=(0,0,-1,0)$
and ${N}_5=(0,0,0,-1)$) do not satisfy the  above conditions.



\vskip 1 cm
\textbf{Affiliations.}\\
National Research University Higher School of Economics,\\
Bolshoi Trekhsvjatitelskii per. 3, Moscow, 109028, Russia\\
e-mails: parus-a@mail.ru, aparusnikova@hse.ru

\end{document}